\documentclass[runningheads,a4paper]{llncs}

\usepackage{amssymb}
\usepackage{graphicx}

\usepackage{amsmath}
\usepackage{url}
\urldef{\mailsa}\path|{d.andreagiovanni,krolikowski,pulaj}@zib.de|
\newcommand{\keywords}[1]{\par\addvspace\baselineskip
\noindent\keywordname\enspace\ignorespaces#1}

\begin{document}

\mainmatter

\title
{
A hybrid primal heuristic for\\Robust Multiperiod Network Design\thanks{
This is the authors' final version of the paper published in: Esparcia-Alc\'azar A., Mora A. (eds), EvoApplications 2014: Applications of Evolutionary Computation, LNCS 8602, pp. 15-26, 2017.
DOI: 10.1007/978-3-662-45523-4\_2.
The final publication is available at Springer via http://dx.doi.org/10.1007/978-3-662-45523-4\_2
}
}

%\title{
%A hybrid primal heuristic for\\Robust Multiperiod Network Design\thanks{This work was partially supported by the \emph{German Research Foundation} (DFG), project \emph{Multiperiod Network Optimization}, by the DFG Research Center Matheon (www.matheon.de), Project B3, and by the \emph{German Federal Ministry of Education and Research} (BMBF), project \emph{ROBUKOM} \cite{BaEtAl13}, grant 05M10PAA.\\
%$^{\dag}$
%}$^{\dag}$
%}

\titlerunning{A hybrid primal heuristic for Robust Multiperiod Network Design}

\author{
Fabio D'Andreagiovanni$^{1,2}$
\and
Jonatan Krolikowski$^{2}$
\and
Jonad Pulaj$^{2}$
}

\authorrunning{F. D'Andreagiovanni, J. Krolikowski, J. Pulaj}

\institute{
$^{1}$DFG Research Center MATHEON, Technical University Berlin\\
Stra{\ss}e des 17. Juni 135, 10623 Berlin\\
\
\\
$^{2}$Department of Optimization, Zuse-Institute Berlin (ZIB)\\
Takustr. 7, 14195 Berlin, Germany\\
\mailsa\\
}

\toctitle{Lecture Notes in Computer Science}
\tocauthor{Authors' Instructions}
\maketitle

\begin{abstract}
We investigate the Robust Multiperiod Network Design Problem, a generalization of the classical Capacitated Network Design Problem that additionally considers multiple design periods and provides solutions protected against traffic uncertainty.
Given the intrinsic difficulty of the problem, which proves challenging even for state-of-the art commercial solvers,
we propose a hybrid primal heuristic based on the combination of ant colony optimization and an exact large neighborhood search. Computational experiments on a set of realistic instances from the SNDlib show that our heuristic can find solutions of extremely good quality with low optimality gap.
\end{abstract}

\section{Introduction}

The design of a telecommunication network can be essentially described as the task of establishing the topology of the network and the technological features (e.g., transmission capacity and rate) of its elements, namely nodes and links.
The dramatic growth that telecommunications have experienced over the last ten years has greatly increased the complexity and difficulty of the corresponding design problems. The growing need for taking into account data uncertainty, such as that of traffic volumes, has made things even more complicated. In this context, the traditional design approach of professionals, based on a combination of trial-and-error and simulation, may lead to arbitrarily bad design solutions and thus the need for optimization-oriented approaches has arisen.

In this paper, we focus on the development of a new Robust Optimization model to tackle traffic uncertainty in a Multiperiod Network Design Problem (MP-NDP). This problem constitutes a natural extension of a classical network design problem, in which we want to decide how to install capacity modules in the network in order to route traffic flows of communications generated by users. The extension implies the design over a time horizon made up of multiple periods. Moreover, traffic uncertainty is taken into account to protect design solutions against deviations of the traffic input data, that may compromise feasibility and optimality of solutions.
To the best of our knowledge, the (MP-NDP) has received little attention and just a few works have investigated it (primarily, \cite{LaNaGe07} and \cite{GaRa12}). These works point out the difficulty of solving multiperiod problems already for just two periods and (easier) splittable-flow routing \cite{LaNaGe07}, and for a pure routing problem in satellite communications \cite{GaRa12}.
Our direct and more recent computational experience confirmed this behaviour, even for instances of moderate size considering a low number of time periods and solved by a state-of-the-art commercial mixed-integer programming solver.

\medskip
\noindent
In this work, our main original contributions are:
\begin{enumerate}
  \item the first Robust Optimization model for Multiperiod Network Design. The formulation is developed to tackle traffic uncertainty, modeling data uncertainty by Multiband Robustness \cite{BuDA12a,BuDA12b,BuDA13,BuDA14}, a new model for Robust Optimization recently introduced to refine the classical Bertsimas-Sim model \cite{BeSi04};
  \item a hybrid solution algorithm, based on the combination of an exact large
    neighborhood search called RINS \cite{DaRoLP05} with ant colony optimization
    \cite{DoDCGa99};
  \item computational experiments over a set of realistic instances derived from SNDlib, the Survivable Network Design Library \cite{SNDlib}, showing that our hybrid algorithm is able to produce solutions of extremely high quality associated with very small optimality gap.
\end{enumerate}
The remainder of this paper is organized as follows: in Section 2, we review a canonical network model for joint routing and capacity installation; in Section 3, we introduce the new formulation for Robust Multiperiod Network Design; in Sections 4 and 5, we present our hybrid metaheuristic and computational results.

\section{Capacitated Network Design}

The \emph{Capacitated Network Design Problem} (CNDP) can be described as follows: given a network and a set of demands whose flows must be routed between vertices of the network, we want to install capacities on network edges and route the flows through the network, so that the capacity constraint of each edge is respected and the total cost of installing capacity is minimized.
The CNDP has been a central and highly studied problem in Network Optimization, that appears in a wide variety of real-world applications. For an exhaustive introduction to the topic, we refer the reader to the well-known book \cite{AhMaOr93}.

The CNDP is commonly formalized in the following way: we are given 1) a network represented by a graph $G(V,E)$, where $V$ is the set of vertices and $E$ the set of edges; 2) a set of commodities $C$, each associated with a traffic flow $d_{c}$ to route from an origin $s_{c}$ to a destination $t_{c}$; 3) a set of admissible paths $P_c$ for routing the flow of each commodity $c$ from $s_{c}$ to $t_{c}$; 4) a cost $\gamma_e$ for installing one module of capacity $\phi > 0$ on edge $e \in E$. Using this notation, we can model the problem as an \emph{integer linear program}:
{
\small
\begin{align}
\min & \sum_{e \in E} \gamma_e \hspace{0.1cm} y_e
&&\mbox{(CNDP-IP)}
\nonumber
\\
&
\sum_{c \in C} \sum_{p \in P_c: \hspace{0.05cm} e \in p} d_{c} \hspace{0.1cm} x_{cp} \leq \phi \hspace{0.1cm} y_e
&&
e \in E
\label{CNDP_cnstr_capacity}
\\
&
\sum_{p \in P_c} x_{cp} = 1
&&
c \in C
\label{CNDP_cnstr_singlePath}
\\
&x_{cp} \in \{0,1\}
&& c \in C, p \in P_c
\nonumber
%\label{CNDP_cnstr_flowVar}
\\
&y_e \in \mathbb{Z}_{+}
&& e \in E \; .
\nonumber
%\label{CNDP_cnstr_capacityVar}
\end{align}
}

\noindent
The problem uses two families of variables: the binary variables $x_{cp}$ (\emph{path-assignment variables}) and the non-negative integer variables $y_e$ (\emph{capacity variables}). A path-assignment variable $x_{cp}$ is equal to $1$ if the entire flow of a commodity $c \in C$ is routed through path $p \in P_c$ and $0$ otherwise. A capacity variable $y_e$ represents instead the number of capacity modules installed on edge $e \in E$.
The objective function minimizes the total cost of installation. Capacity constraints \eqref{CNDP_cnstr_capacity} impose that the summation of all flows routed through an edge $e \in E$ must not exceed the capacity installed on $e$ (equal to the number of installed modules represented by $y_e$ multiplied by the capacity $\phi$ granted by a single module). Constraints \eqref{CNDP_cnstr_singlePath} impose that flow of each commodity $c \in C$ must be routed through a single path.

\smallskip
\noindent
\textbf{Remark 1.}
This is an unsplittable version of the CNDP, namely the traffic flow of a commodity $c \in C$ \emph{cannot be split} over multiple paths going from $s_{c}$ to $t_{c}$, but must be routed on exactly one path.
Moreover, the set of feasible routing paths $P_c$  of each commodity is pre-established and constitutes an input of the problem.
This is in line with other works based on industrial cooperations (e.g., \cite{BlGrWe07}) and with our experience \cite{BaEtAl13}, in which a network operator typically considers just a few paths that meet its own specific business and quality-of-service considerations.

\section{Multiband-Robust Multiperiod Network Design}

We introduce now a generalization of the CNDP, designing the network over multiple time periods and taking into account  traffic uncertainty.
The multiperiod design requires the introduction of a time horizon made up of a set of elementary time periods $T = \{1,2, \ldots, |T|\}$. From a modeling point of view, in the optimization problem we simply need to add a new index $t \in T$ to the decision variables, to represent routing and capacity installation decisions taken in each period (we stress however that this greatly increases the size and complexity of the problem).

Concerning traffic uncertainty,
we assume that for each commodity $c \in C$ the demand $d_c$ is \emph{uncertain}, i.e. its value is not known exactly, but lies in a known range. More specifically,  we assume to know a nominal value of traffic $\bar{d}_{c}$ and maximum negative and positive deviations $\delta_{c}^{-}, \delta_{c}^{+}$ from it. The actual value $d_c$ thus belongs to the interval:
$
\hspace{0.1cm}
d_c
\hspace{0.1cm} \in \hspace{0.1cm}
[\bar{d}_{c} - \delta_{c}^{-}, \hspace{0.2cm} \bar{d}_{c} + \delta_{c}^{+}]
$.

\medskip
\noindent
\textbf{Example 1 (traffic uncertainty).}
We are given two commodities $c_1, c_2$ with nominal traffic demands $\bar{d}_{c_1} = 100$ Mb, $\bar{d}_{c_2} = 150$ Mb and we know that these values may deviate up to 10\%. So the maximum negative and positive deviations for $c_1, c_2$ are $\delta_{c_1}^{-} = \delta_{c_1}^{+} = 10$ Mb, $\delta_{c_2}^{-} = \delta_{c_2}^{+} = 15$ Mb, respectively. The actual values of traffic are then $d_{c_1} \in [90,110]$ Mb, $d_{c_2} \in [135,165]$ Mb.

\medskip
\noindent
The presence of uncertain data in an optimization problem can be very tricky: it is well-known that even small variations in the value of input data may make an optimal solution heavily suboptimal, whereas feasible solutions may reveal to be infeasible and thus completely useless in practice \cite{BeElNe09}. As a consequence, in our case we cannot optimize just using the nominal demand values $\bar{d}_{c}$, but we must take into account the possibility that demands will vary in the ranges $[\bar{d}_{c} - \delta_{c}^{-}, \hspace{0.2cm} \bar{d}_{c} + \delta_{c}^{+}]$ that we have characterized.
We illustrate the bad effects of input data deviations by providing an example.

\medskip
\noindent
\textbf{Example 2 (infeasibility caused by deviations).}
Consider again the commodities of Example 1 and suppose that in some link we have installed exactly the capacity to handle the sum of their nominal values (i.e., we have installed $100+150$ Mb of capacity).
This capacity dimensioning neglects that the demands may deviate up to 10\%. It is sufficient that one demand increases, while the other remains the same to violate the capacity constraint of the link, making the design solution infeasible in practice.

\medskip
\noindent
Over the years, many methods such as \emph{Stochastic Programming} and \emph{Robust Optimization} have been proposed in literature for dealing with data uncertainty in optimization problems. We refer the reader to \cite{BeElNe09} for a general discussion about data uncertainty and its effects and for an overview of the most studied methodologies to deal with them.

In this paper, we tackle data uncertainty by Robust Optimization (RO), a methodology that has gained a lot of attention over the last decade \cite{BeElNe09,BeSi04}. RO essentially takes into account data uncertainty by including additional hard constraints in the optimization problem. These constraints eliminate those solutions that are not protected against deviations of the input data from their nominal values. So a robust optimization problem considers only those solutions that are completely protected against specified data deviations. The data deviations that are considered are specified through a so-called \emph{uncertainty set}.
More formally, suppose that we are given a generic linear program:
%\\
%\centerline{
    $$v = \max
    \hspace{0.1cm}
        c' x
        \hspace{0.1cm}
        \mbox{ with }
        \hspace{0.05cm}
        x \in {\cal F} =
        \{
        A \hspace{0.05cm} x \leq b,
        \hspace{0.15cm}
        x \geq 0
        \}$$
%}
%
and that the coefficient matrix $A$ is uncertain, i.e. we do not know exactly the value of its entries. However, we are able to identify a family $\cal A$ of coefficient matrices that represent possible valorizations of the uncertain matrix $A$, i.e. $A \in \cal A$. This family represents the uncertainty set of the robust problem. Then we can produce a \emph{robust optimal solution}, i.e. a solution that is protected against data deviations, by considering the \emph{robust counterpart} of the original problem:
%
%\\
%\centerline{
$$
v^{{\cal R}} = \max
\hspace{0.1cm}
    c' x
    \hspace{0.1cm}
    \mbox{ with }
    \hspace{0.05cm}
    x \in {\cal R} =
    \{
    \tilde{A} \hspace{0.1cm} x \leq b
    \hspace{0.3cm}
    \forall \tilde{A} \in {\cal A},
    \hspace{0.15cm}
    x \geq 0
    \} \; .
$$
The feasible set ${\cal R}$ of the robust counterpart contains only those solutions  that are feasible for all the coefficient matrices in the uncertainty set ${\cal A}$. Therefore, ${\cal R}$ is a subset of the feasible set of the original problem, i.e. ${\cal R} \subseteq {\cal F}$. The choice of the coefficient matrices included in ${\cal A}$ should reflect the risk aversion of the decision maker.

Providing protection entails the so-called \emph{price of robustness}, namely a deterioration of the optimal value of the robust counterpart w.r.t. the optimal value of the original problems (i.e., $v^{{\cal R}} \leq v$). This is a consequence of restricting the feasible set to only robust solutions. The price of robustness reflects the features of the uncertainty set: uncertainty sets expressing higher risk aversion will take into account more severe and unlikely deviations, leading to higher protection but also higher price of robustness; conversely, uncertainty sets expressing risky attitudes will tend to neglect improbable deviations, offering less protection but also a reduced price of robustness.

\medskip
\noindent
\textbf{Example 3 (protection against deviations).}
Following example 2, a simple way to grant protection would be to install sufficient capacity to deal with the peak deviations of each commodity. So we should install 110+165 Mb of capacity.

\subsection{A Robust Optimization model for traffic-uncertain Multiperiod Network Design}

If we denote by ${\cal D}$ the uncertainty set associated with the demands of the commodities, we can finally state the general form of the robust counterpart of the multiperiod network design problem as follows:
{
\small
\begin{align}
\min & \sum_{e \in E} \sum_{t \in T} \gamma_{et} \hspace{0.1cm} y_{et}
&&
%\mbox{(CNDP)}
\nonumber
\\
&
\sum_{c \in C} \sum_{p \in P_c: \hspace{0.05cm} e \in p} \bar{d}_{ct} \hspace{0.1cm} x_{cpt}
+ DEV_{et}(x, {\cal D})
\leq \phi \sum_{\tau = 1}^{t} y_{e\tau}
&&
e \in E, t \in T
\label{Rob-CNDP_robustcnstr_capacity}
\\
&
\sum_{p \in P_c} x_{cpt} = 1
&&
c \in C, t \in T
\nonumber
%\label{CNDP_cnstr_singlePath}
\\
&x_{cpt} \in \{0,1\}
&& c \in C, p \in P_c, t \in T
\nonumber
%\label{CNDP_cnstr_flowVar}
\\
&y_{et} \in \mathbb{Z}_{+}
&& e \in E, t \in T \; .
\nonumber
%\label{CNDP_cnstr_capacityVar}
\end{align}
}

\noindent
Besides the addition of a new index $t \in T$ in the decision variables to represent decisions taken in each time period, the modifications of the model concentrates in the robust capacity constraints \eqref{Rob-CNDP_robustcnstr_capacity}. Each of these constraints considers:
1) the sum of \emph{nominal} traffic demands $\bar{d}_{ct}$ of commodities using the edge $e$ in period $t$;
2) the overall maximum positive deviation $DEV_{et}(x, {\cal D})$ that demands may experience on edge $e$ in period $t$ and are allowed by the uncertainty set ${\cal D}$ for a routing vector $x$;
3) the overall capacity installed in $e$ since the first period of the horizon (so we sum up the integer variables $y_{e\tau}$ from period $1$ to $t$ and multiply them by the basic capacity $\phi$ of a module).

\smallskip
\noindent
\textbf{Structuring the uncertainty set ${\cal D}$.}
We now have a general definition of the robust counterpart of the multiperiod problem. A question that is still open is how to structure the uncertainty set ${\cal D}$ and deciding which deviations from the nominal traffic values $\bar{d}_{ct}$ to take into account to produce robust solutions.

To characterize ${\cal D}$, we use \emph{Multiband Robustness}, a new model for Robust Optimization recently introduced to refine and generalize the classical $\Gamma$-robustness model by Bertsimas and Sim \cite{BeSi04}, while maintaining its accessibility and tractability. For a detailed explanation of Multiband Robustness we refer the reader to \cite{BuDA12a,BuDA12b,BuDA13,BuDA14}. Here we directly discuss the adaption of the model to our specific case.

\smallskip
\noindent
According to the multiband framework, we build a \emph{multiband uncertainty set} as follows:
\begin{enumerate}
  \item for each commodity $c \in C$ and time period $t \in T$, we know the nominal value $\bar{d}_{ct}$ of the traffic coefficient and maximum negative and positive deviations $\delta_{ct}^{-},\delta_{ct}^{+}$ from it. The actual value $d_{ct}$ is then such that $d_{ct} \in [\bar{d}_{ct} - \delta_{ct}^{-}, \hspace{0.1cm} \bar{d}_{ct} + \delta_{ct}^{+}]$;
  \item the overall deviation range $[\bar{d}_{ct} - \delta_{ct}^{-}, \hspace{0.1cm} \bar{d}_{ct} + \delta_{ct}^{+}]$ of each coefficient $d_{ct}^{}$ is partitioned into $K$ bands, defined on the basis of $K$ deviation values:
      \\
        $
        -\infty<
        - {\delta_{ct}^{-} = \delta_{ct}^{K^{-}}<\cdots<\delta_{ct}^{-1}
        \hspace{0.1cm}<\hspace{0.1cm}\delta_{ct}^{0}=0\hspace{0.1cm}<\hspace{0.1cm}
        \delta_{ct}^{1}<\cdots<\delta_{ct}^{K^{+}}} = \delta_{ct}^{+}
        <+\infty \; ;
        $
  \item through these deviation values, $K$ deviation bands are defined, namely:
    a set of positive deviation bands $k\in \{1,\ldots,K^{+}\}$ and  a set of negative deviation bands $k \in \{K^{-}+1,\ldots,-1,0\}$, such that a band $k\in \{K^{-}+1,\ldots,K^{+}\}$ corresponds to the range $(d_{ct}^{k-1},d_{ct}^{k}]$, and band $k = K^{-}$
    corresponds to the single value $d_{ct}^{K^{-}}$;
  \item for each capacity constraint $\eqref{Rob-CNDP_robustcnstr_capacity}$ defined for an edge $e \in E$, period $t \in T $ and band $k\in K$, a value $\theta_{etk} \geq 0$ is introduced to represent the number of traffic coefficients of the constraint whose value deviates in band $k$. Of course, $\theta_{etk} \geq 0$ must be less or equal than the number of traffic coefficients that are present in the constraint.
\end{enumerate}

\noindent
Given the previous characterization of the multiband uncertainty set, the maximum positive deviation of traffic $DEV_{et}(x, {\cal D})$ of a constraint \eqref{Rob-CNDP_robustcnstr_capacity} can be found by solving a binary linear program (see \cite{BuDA12a} for details). Since the polytope associated with the binary program is shown to be integral, by considering its relaxation and by exploiting strong duality, it is possible to reformulate the original trivial robust counterpart as the following \emph{linear and compact} robust counterpart (we refer the reader to \cite{BuDA12a} for a formal proof of the result):
{
\small
\begin{align}
\min & \sum_{e \in E} \sum_{t \in T} \gamma_{et} \hspace{0.1cm} y_{et}
&&
\mbox{(Rob-MP-CNDP)}
\nonumber
\\
&
\sum_{c \in C} \sum_{p \in P_c: \hspace{0.05cm} e \in p} \bar{d}_{ct} \hspace{0.1cm} x_{cpt}
\hspace{0.1cm} +
&&
\nonumber
\\
&
+
%\hspace{0.1cm}
\sum_{k \in K} \theta_{etk} \hspace{0.05cm} w_{etk}
+ \sum_{c \in C} \sum_{p \in P_c: \hspace{0.05cm} e \in p} z_{ecpt}
\leq
\phi \sum_{\tau = 1}^{t} y_{e\tau}
&&
e \in E, t \in T
\nonumber
%\label{COMPACTRobCNDP_cnstr_capacity}
\\
&
w_{etk} + z_{ecpt} \geq \delta_{ctk} \hspace{0.1cm} x_{cpt}
&&
e \in E, c \in C, p \in P_c: e \in p,
\nonumber
\\
&
&&
t \in T, k \in K
\label{RobCNDP_cnstr_dual}
%\\
%&
%w_{etk} + z_{ecpt} \geq \delta_{ctk} \hspace{0.1cm} x_{cpt}
%&&
%e \in E, c \in C, p \in P_c: e \in p, t \in T, k \in K
%\label{RobCNDP_cnstr_dual}
\\
&w_{etk} \in \mathbb{R}
&& e \in E, t \in T, k \in K
\label{RobCNDP_dual1}
\\
&z_{ecpt} \geq 0
&& e \in E, c \in C, p \in P_c: e \in p, t \in T
\label{RobCNDP_dual2}
%\end{align}
%
%\begin{align}
\\
&
\sum_{p \in P_c} x_{cpt} = 1
&&
c \in C, t \in T
\nonumber
%\label{CNDP_cnstr_singlePath}
\\
&x_{cpt} \in \{0,1\}
&& c \in C, p \in P_c, t \in T
\nonumber
%\label{CNDP_cnstr_flowVar}
\\
&y_{et} \in \mathbb{Z}_{+}
&& e \in E, t \in T \; ,
\nonumber
%\label{CNDP_cnstr_capacityVar}
\end{align}
}

\noindent
This formulation includes additional constraints \eqref{RobCNDP_cnstr_dual} and variables \eqref{RobCNDP_dual1},\eqref{RobCNDP_dual2} which are derived from the dualization operation that allow to linearly reformulate the original (non-linear) problem including the term $DEV_{et}(x, {\cal D})$ in each capacity constraint (see \cite{BuDA12a} for details).

In principle, we can get a robust optimal solution for (Rob-MP-CNDP) by using any commercial mixed-integer programming software. However, as showed in the computational experiment section, getting feasible solutions to this problem may be a challenge even for a state-of-the-art solver like IBM ILOG CPLEX (http://www-01.ibm.com). In the next section, we thus propose a hybrid exact-ant colony primal heuristic that is able to find solutions of very high quality.

\section{A hybrid primal heuristic for the Rob-MP-CNDP}
\label{sec:ACO}

Attracted by the effectiveness of MIP-based and bio-inspired heuristics in hard network design problems (see, for example \cite{DA11,DAMaSa11,DoDCGa99,Ka13}),
we present an original hybrid primal heuristic based on the combination of Ant Colony Optimization (ACO) and an exact large neighbourhood search.
ACO is a metaheuristic originally proposed by Dorigo and colleagues
for combinatorial optimization \cite{DoMaCo96}
and later extended to integer and continuous problems (e.g., \cite{DoDCGa99}). Over the years several refinements of the basic algorithm have been proposed (e.g., \cite{GaMoWe12,Ma99}).
ACO was inspired by the behaviour of ants searching for food and is essentially based on the definition of a cycle where a number of feasible solutions are iteratively built in parallel, using information about solutions built in previous executions of the cycle. An ACO algorithm presents the following general structure:
{
\small
\begin{enumerate}
  \item UNTIL an arrest condition is reached DO
      \hspace{2.0cm} (Gen-ACO)
      \begin{enumerate}
      \item Ant-based solution construction
      \item Pheromone trail update
    \end{enumerate}
    \item Daemon actions
\end{enumerate}
}
\noindent
We now proceed to detail each phase of the previous sketch for our hybrid ACO-exact algorithm for the (Rob-MP-CNDP). Our approach is hybrid since the canonical ACO construction phase is followed by a daemon-action phase, based on an exact large neighborhood search formulated as a mixed-integer linear program.

\smallskip
\noindent
\textbf{Ant-based solution construction.}
In the step 1 of the cycle, $m \geq 0$ \emph{ants} are defined and each ant iteratively builds a feasible solution for the optimization problem. At every iteration, the ant is in a \emph{state} corresponding with a \emph{partial solution} of the problem and can further complete the solution by making a \emph{move} and thus fixing the value of a new non-fixed variable. The move is chosen probabilistically, evaluating pheromone trail values. For a more detailed description of the elements and actions of step 1, we refer the reader to the paper \cite{Ma99} by Maniezzo. This paper presents ANTS, an improved ANT algorithm that we have taken as reference for our work. We considered ANTS particularly attractive as it proposes a series of improvements for ACO that allow to better exploit polyhedral information about the problem. Furthermore, ANTS is based on a reduced number of parameters and uses more efficient mathematical operations.

Before describing how our ANTS implementation is structured, we make some preliminary considerations. The formulation (Rob-MP-CNDP) is based on four families of variables: 1) the path assignment variables $x_{cpt}$; 2) the capacity variables $y_{et}$; 3-4) the auxiliary variables $w_{etk}$, $z_{ecpt}$ coming from robust dualization.
Though we have to deal with four families, we can notice that routing decisions taken over the time horizon entirely determine the capacity installation of minimum cost.
Indeed, once the values of all path assignment variables are fixed, the routing is completely established and the worst traffic deviation term $DEV_{et}(x, {\cal D})$ can be efficiently derived without the auxiliary variables $w_{etk}$, $z_{ecpt}$ \cite{BuDA12a,BuDA13}. So we can derive the total traffic $D_{et}$ sent over an edge $e$ in period $t$ in the worst case. The minimum cost installation can then be derived through a sequential evaluation from period 1 to period $|T|$, keeping in mind that we must have $\left\lceil \frac{D_{et}}{\phi} \right\rceil$ capacity modules on $e$ in $t$ to accommodate the traffic.
As a consequence, in the ant-construction phase we can limit our attention to the binary assignment variables and we introduce the concept of \emph{routing state}.
\begin{definition}
  Routing state (RS): let $P = \bigcup_{c \in C} P_c$ and let $R \subseteq C \times P \times T$ be the subset of triples $(c,p,t)$ representing the assignment of path $p \in P_c$ to commodity $c$ in period $t \in T$. A \emph{routing state} is an assignment of paths to a subset of commodities in a subset of time periods which excludes that multiple paths are assigned to a single commodity. Formally:
%  \\
%    \centerline{
        $$
        RS \subseteq R:
       \hspace{0.2cm}
       \not \exists (c_1,p_1,t_1), (c_2,p_2,t_2) \in RS:
       \hspace{0.1cm}
       c_1 = c_2
       \hspace{0.1cm} \wedge \hspace{0.1cm}
       p_1,p_2 \in P_{c_1}
     \hspace{0.1cm} \wedge \hspace{0.1cm}
       t_1 = t_2 \; .
       $$
%    }
\end{definition}

\noindent
We say that a routing state RS is \emph{complete} when it specifies the path used by each commodity in each time period (thus $|RS| = |C||T|$). Otherwise the RS is called \emph{partial} and we have $|RS| < |C||T|$).
\\
In the ANTS algorithm that we propose, we decided to assign paths considering time periods and commodities in a pre-established order. Specifically, we establish the routing in each time period separately, starting from $t=1$ and continuing up to $t=|T|$, and in each time period commodities are sorted in descending order w.r.t. their nominal traffic demand. Formally, this can be sketched through the following cycle that builds a \emph{complete routing state}:
{
\small
\center
\begin{description}
\item FOR $t := 1$ TO $|T|$ DO
      \begin{enumerate}
        \item sort $c \in C$ in descending order of $\bar{d}_{ct}$.
        \item FOR (sorted $c \in C$) DO
          \begin{enumerate}
            \item assign a single path $p \in P_c$ to $c$;
          \end{enumerate}
        END FOR
      \end{enumerate}
\item END FOR
\end{description}
}

\medskip
\noindent
For an iteration $(t,c)$ of the above nested cycles, the assignment of a path to a commodity corresponds with an ant moving from a partial routing state $RS_i$ to a partial routing state $RS_j$ such that:
$
RS_j = RS_i \cup \{(c,p,t)\} \hspace{0.1cm} \mbox{ with } p \in P_c\;  .
$
\\
We note that by the definition of routing state a sequence of moves is actually a sequence of fixings of decision variables, as done in \cite{Ma99}.

The probability that an ant $k$ moves from a routing state $i$ to a more complete routing state $j$, chosen among a set of feasible routing states, is defined by the improved formula of \cite{Ma99}:
$$
p_{ij}^{k} = \frac{\alpha \hspace{0.1cm} \tau_{ij} + (1-\alpha) \hspace{0.1cm} \eta_{ij}}
                {\sum_{f \in F} \alpha \hspace{0.1cm} \tau_{if} + (1-\alpha) \hspace{0.1cm} \eta_{if}}
\; ,
$$
where $\alpha \in [0,1]$ is a parameter assessing the relative importance of trail and attractiveness.
As discussed in \cite{Ma99}, the trail values $\tau_{ij}$  and the attractiveness values $\eta_{ij}$ should be provided by suitable lower bounds of the considered optimization problem. In our particular case:
1) $\tau_{ij}$ is derived from the values of the variables in the solution associated with the linear relaxation of the robust counterpart (Rob-MP-CNDP);
2) $\eta_{ij}$ is equal to the optimal solution of the linear relaxation
of the nominal multiperiod network design problem, i.e. the problem that does not consider the traffic uncertainty. The optimum of this problem can be quickly computed and its computation becomes faster as more variables are fixed.

\smallskip
\noindent
\textbf{Daemon actions: Relaxation Induced Neighborhood Search.}
At the end of the ant-construction phase, we try to improve the quality of the feasible solution found by executing an \emph{exact local search} in a \emph{large neighborhood}. In particular, we adopt a modified \emph{Relaxation Induced Neighborhood Search} (RINS) (see \cite{DaRoLP05} for an exhaustive description of the method).
Let $(\bar{x},\bar{y})$ be a feasible solution of (Rob-MP-CNDP) found by an ant and $(x^{\small LR},y^{\small LR})$ be an optimal (continuous) solution of the linear relaxation of (Rob-MP-CNDP)
Moreover, let $(\bar{x},\bar{y})_j, (x^{\small LR},y^{\small LR})_j$ denote the $j$-th component of the vectors. Our modified RINS \emph{(mod-RINS)} solves a sub-problem of (Rob-MP-CNDP) where:
\begin{enumerate}
    \item we fix the variables $x$ whose value in $(\bar{x},\bar{y})$ and $(x^{\small LR},y^{\small LR})$ differs of at most $\epsilon > 0$, i.e.:
        \begin{description}
          \item $\bar{x}_j = 0$
          $\hspace{0.1cm} \cap \hspace{0.1cm}$
          $x^{\small LR}_j \leq \epsilon$
          $\hspace{0.2cm} \Longrightarrow \hspace{0.2cm}$
          $x_j = 0$
          \item $\bar{x}_j = 1$
          $\hspace{0.1cm} \cap \hspace{0.1cm}$
          $x^{\small LR}_j \geq 1 - \epsilon$
          $\hspace{0.2cm} \Longrightarrow \hspace{0.2cm}$
          $x_j = 1$
        \end{description}
      %\item set an objective cutoff based on the value of the ant solution $(\bar{x},\bar{y})$;
      \item impose a solution time limit of $T$.
\end{enumerate}
A time limit is imposed since the subproblem may be difficult to solve, so the exploration of the neighbourhood may need to be truncated. Note that in point 1 we generalize the fixing rule of RINS, in which $\epsilon = 0$.

\smallskip
\noindent
\textbf{Pheromone trail update.}
At the end of each ant-construction phase $h$, the pheromone trails of a move $\tau_{ij}(h-1)$ are updated according to an improved formula proposed in \cite{Ma99}:
\begin{equation}
\small
\label{pheroFormula}
\tau_{ij}(h) \hspace{0.05cm} = \hspace{0.05cm} \tau_{ij}(h-1)
\hspace{0.05cm} + \hspace{0.05cm} \sum_{k=1}^{m} \tau_{ij}^k
\hspace{0.5cm} \mbox{ with } \hspace{0.1cm} \tau_{ij}^k =
\hspace{0.1cm}
\tau_{ij}(0)
\cdot \left(
1 - \frac{z_{curr}^k - LB}{\bar{z} - LB}
\right) ,
\end{equation}

\noindent
where the values $\tau_{ij}(0)$ and $LB$ are set by using the linear relaxation of (Rob-MP-CNDP): $\tau_{ij}(0)$ is set equal to the values of the corresponding optimal decision variables and $LB$ equal to the optimal value of the relaxation. Additionally, $z_{curr}^k$ is the value of the solution built by ant $k$ and $\bar{z}$ is the moving average of the values of the last $\psi$ feasible solutions built.
As noticed in \cite{Ma99}, adopting formula \eqref{pheroFormula} allows to replace the pheromone evaporation factor, a tricky parameter, with the moving average $\psi$ whose setting has been shown to be much less critical.

\medskip
\noindent
Algorithm 1 details the structure of our original hybrid exact-ACO algorithm. The algorithm includes an outer loop repeated until a time limit is reached. At each execution of the loop, an inner loop defines $m$ ants to build the solutions. Pheromone trail updates are done at the end of each execution of the inner loop. Once the ant construction phase is over, mod-RINS is applied so to try to get an improvement.

\medskip
\noindent
{
\small
\noindent
\textbf{Algorithm 1. Hybrid ACO-exact algorithm for (Rob-MP-CNDP).}
\begin{enumerate}
  \item Compute the linear relaxation of (Rob-MP-CNDP) and initialize
        the values of $\tau_{ij}(0)$ by it.
  \item UNTIL time limit is reached DO
      \begin{enumerate}
        \item FOR $\mu := 1$ TO $m$ DO
          \begin{enumerate}
            %\item Compute the linear relaxation of (strongBM-SPCAP) and use it to initialize the values $\eta_{ij}$;
            \item build a complete routing state;
                %\ref{subsec:antConstruction};
                %
                %(by values $\eta_{ij}$ computed through (strongBM-SPCAP));
            \item derive a complete feasible solution for (Rob-MP-CNDP);
          \end{enumerate}
        END FOR
        \item Update $\tau_{ij}(t)$ according to (\ref{pheroFormula}).
      \end{enumerate}
      \item apply mod-RINS to the best feasible solution.
\end{enumerate}
}

\section{Experimental results}  \label{sec:computations}

We tested the performance of our hybrid algorithm on a set of 15 instances based on realistic network topologies from the SNDlib \cite{SNDlib} defined in collaboration with industrial partners from former and ongoing projects. The experiments were performed on a machine with a 2.40 GHz quad-core processor and 16 GB of RAM and using IBM ILOG CPLEX 12.4. All the instances led to very large and hard to solve robust multiperiod network design problems. We observed that even a state-of-the-art solver like CPLEX had troubles identifying good feasible solutions and in all the cases the final optimality gap was over 90\%. In contrast, as clear from Table \ref{table:results}, in most cases our hybrid primal heuristic was able to find very high quality solutions associated with very low optimality gaps.

After executing preliminary tests, we found that an effective setting of the parameters of the heuristic was: $\alpha = 0.5$ (balancing attractiveness and trail level), $m = 3$ ants, $\psi = m$ (width of the moving average equal to the number of ants), $\epsilon = 0.1$ (tolerance of fixing in mod-RINS), $T = 20$ minutes (time limit in mod-RINS). The overall time limit for the execution of the heuristic was 5 hours. The same time limit was imposed on CPLEX when used to solve the robust counterpart (Rob-MP-CNDP). Each commodity admits 5 feasible paths, i.e. $|P_c| = 5, \forall c \in C$ and 3 positive deviations bands including the null deviation band. For each instance, in Table \ref{table:results} we report its ID and features ($|V|$ = no. vertices, $|E|$ = no. edges, $|C|$ = no. commodities, $|T|$ = no. time periods). Moreover, we show  the performance of the hybrid solution approach, that is denoted by the three measures $c^{*}$(ACO), $c^{*}$(ACO+RINS), gapAR\%, which  represent the value of the best solution found by pure ACO, the value of the best solution found by ACO followed by RINS and the corresponding final optimality gap). We also show the performance of CPLEX, which is denoted by measures $c^{*}$(IP) and gapIP\% representing the value of the best solution found and the corresponding final optimality gap.

The best solutions found by our hybrid algorithm have in most cases a value that is at least one order of magnitude better than those found by CPLEX (2700\% better on average). The results are of very high quality and, given the very low optimality gap, we can suppose that some of these solutions are actually optimal. We notice that in most cases executing RINS after the ant-construction phase can remarkably improve the value of the best solution found by the ants.
\begin{table}
\label{table:results}
\caption{Experimental results}
\small
\begin{center}
\begin{tabular}{c c c c @{\quad} c @{\quad}
 c c c c @{\quad} c}
\hline
\hline
ID & $|V|$ & $|E|$ & $|C|$ & $|T|$ & $c^{*}$(ACO) & $c^{*}$(ACO+RINS) & gapAR\% &
$c^{*}$(IP) & gapIP\%
\\ [2pt]
\hline
\hline
&  &    &   &  5 & 1.16E07 & 5.68E06 & 29.8 & 1.37E08 & 97.1
\\
Germany50 & 50 &  88 &  662 &  7 & 2.12E07 & 9.02E6 & $15.5$ & 3.48E08 & $97.8$
\\
&  &    &   &  10 & 6.66E07 & 5.75E08 & 96.2 & 1.25E09 & 98.2
\\
\hline
&  &   &   &  5 & 5.89E06 & 2.34E06 & 1.3 & 9.52E07 & 97.6
\\
Pioro40 & 40 &  89 &  780 & 7 & 1.42E07 & 5.10E06 & 3.1 & 2.40E08 & 97.9
\\
&  &   &   & 10 & 4.78E07 & 1.62E07 & 0.4 & 8.45E08 & 98.1
\\
\hline
&  &   &  & 5 & 6.41E06 & 3.04E06 & 23.0 & 6.01E07 & 96.1
\\
Norway & 27 &  51 &  702 &  7 & 1.44E07 & 5.73E06 & 12.8 & 1.47E08 & 96.6
\\
&  &   &  & 10 & 4.91E07 & 1.74E07 & 7.7 & 5.15E08 & 96.9
\\
\hline
&  &  &  & 5 & 1.55E06 & 6.04E05& 2.2 & 1.74E07& 96.6
\\
Geant & 22 & 36 & 462 & 7 & 3.61E06 & 1.29E06 & 1.6 & 4.32E07 & 97.1
\\
&  &  &  & 10 & 1.23E07 & 4.30E06 & 0.5 & 1.24E08 & 96.5
\\
\hline
&  &   &  & 5 & 2.55E05 & 1.02E05 & 4.9 & 1.50E06 & 93.5
\\
France & 25 & 45 & 300 & 7 & 5.97E05 & 2.18E05 & 2.2 & 3.01E06 & 92.9
\\
&  &  &  & 10 & 2.00E06 & 6.81E05 & 1.0 & 1.62E07 & 95.8
\\
\hline
\hline
\end{tabular}
\end{center}
\end{table}

\section{Conclusion and future work}  \label{sec:end}
We studied a Robust Optimization model for the Multiperiod Network Design Problem to tackle uncertainty of traffic demands. Robust solutions are deterministically protected against deviations of input traffic data, that may compromise the quality of produced solutions.
The increase in complexity and dimension of the problem caused by considering multiple periods and robustness prevents state-of-the-art commercial solvers from finding good quality solutions, so
we have defined a hybrid heuristic based on the combination of ant colony optimization and an exact large neighborhood search. Computational experiments on a set of realistic instances from the SNDlib showed that our heuristic can find solutions of extremely good quality.
As future work, we plan to refine the heuristic (for example, by improving the ant-construction phase) and to integrate it with a branch-and-cut algorithm to enhance its computational performance.

\section*{Acknowledgments}
This work was partially supported by the \emph{German Research Foundation} (DFG), project \emph{Multiperiod Network Optimization}, by the DFG Research Center Matheon (www.matheon.de), Project B3, and by the \emph{German Federal Ministry of Education and Research} (BMBF), project \emph{ROBUKOM} \cite{BaEtAl13}, grant 05M10PAA.

%\bibliographystyle{splncs}
%\bibliography{DAndreagiovanni}

\end{document}